\documentclass[11pt,reqno]{amsart}
\usepackage[toc,page]{appendix}
\usepackage{url}
\usepackage{amsmath} 
\usepackage{amssymb}
\usepackage{amsthm}
\usepackage{bm, bbm}
\usepackage[dvipsnames]{xcolor}
\usepackage{mathtools}
\usepackage{caption}
\usepackage{xfrac}  
\usepackage{enumitem} 
\usepackage[normalem]{ulem}

\usepackage{adjustbox}

\makeatletter
\def\section{\@startsection{section}{1}%
	\z@{.7\linespacing\@plus\linespacing}{.5\linespacing}%
	{\normalfont\bfseries\centering}}
\makeatother

\theoremstyle{plain}
\newtheorem{theorem}{Theorem}[section] 
\newtheorem{corollary}[theorem]{Corollary} 
\newtheorem{lemma}[theorem]{Lemma}
\newtheorem{proposition}[theorem]{Proposition}
\newtheorem{remark}[theorem]{Remark}

\newtheorem{example}[theorem]{Example}

\definecolor{unia-purple}{RGB}{173, 0, 124}
\definecolor{light-gray}{rgb}{0.8889,0.8889,0.8889}
\definecolor{cpl3}{rgb}{1.0,0.0,0.0}
\definecolor{cpl1}{rgb}{0.8889,0.4356,0.2781}
\definecolor{cpl2}{rgb}{0.0,0.6056,0.9787}
\definecolor{cpl4}{rgb}{0.2422,0.6433,0.3044}
\definecolor{cpl6}{rgb}{0.2, 0.2, 0.2}

\usepackage{scalerel,stackengine}
\usepackage[
  backend=biber,
  style=numeric,
  sorting=nty,
  giveninits=true,
  maxnames=5,
  isbn=false
]{biblatex}
\renewbibmacro*{in:}{}

\DeclareFieldFormat[
  article,incollection,inproceedings,thesis,misc,unpublished
]{title}{#1}

\AtEveryBibitem{%
  \ifentrytype{article}
    {\clearfield{number}\clearfield{url}}
    {}%
}

\usepackage{hyperref}
\hypersetup{
	colorlinks=true,
	linkcolor=blue,
	filecolor=blue,      
	urlcolor=blue,
	citecolor=blue
}

\newcommand\R{\mathbb R}

\DeclareMathOperator{\ggrad}{grad}
\DeclareMathOperator{\Hess}{Hess}
\DeclareMathOperator{\Exp}{Exp}
\DeclareMathOperator{\Log}{Log}
\DeclareMathOperator{\D}{\nabla}
\DeclareMathOperator{\Ddual}{\nabla^\ast}

\DeclareMathOperator{\diag}{diag}
\DeclareMathOperator{\dom}{dom}

\newcommand{\VF}[1]{\mathfrak{X}(#1)}
\newcommand{\Dsym}{D^{\mathrm{sym}}}
\newcommand{\avM}{\bar H}
\newcommand{\op}{\mathrm{op}}

\newcommand{\LC}{{\rm LC}}

\newcommand{\vphi}{\varphi}
\newcommand{\la}{\langle}
\newcommand{\ra}{\rangle}
\DeclareMathOperator{\intr}{int}

\makeatletter
\def\author@andify{%
	\nxandlist {\unskip ,\penalty-1 \space\ignorespaces}%
	{\unskip {} \@@and~}%
	{\unskip \penalty-2 \space \@@and~}%
}
\renewcommand{\andify}{%
	\nxandlist{\unskip, }{\unskip{} \@@and~}{\unskip \space \@@and~}}
\makeatother

\definecolor{myBlue}{rgb}{0.0, 0.0, 0.8}
\definecolor{myOrange}{RGB}{225,92,22} 
\definecolor{myGreen2}{RGB}{114,175,30} 
\definecolor{myPurple}{rgb}{0.494,0.184,0.556}%

\newif\ifshowcomments
\showcommentstrue

\ifshowcomments
\newcommand{\YE}[1]{{\color{WildStrawberry} #1}}
\newcommand{\SP}[1]{{\color{myGreen2} #1}}
\newcommand{\todo}[1]{{\color{red}#1}}
\newcommand{\del}[1]{{\color{gray!50}#1}}

\else
\newcommand{\YE}[1]{}
\newcommand{\SP}[1]{}
\newcommand{\todo}[1]{}
\newcommand{\del}[1]{}
\fi

\title[Relative Smoothness]{On Relative Smoothness in Dually Flat Geometry}

\author[Y. Elshiaty]{Yara Elshiaty}
\author[S. Petra]{Stefania Petra}

\thanks{Y. Elshiaty: Institute for Mathematics,
Heidelberg University, Im Neuenheimer Feld 205,
69120 Heidelberg, Germany.
Email: \texttt{elshiaty@math.uni-heidelberg.de}.}

\thanks{S. Petra: Institute for Mathematics, IWR \& MIISM,
Heidelberg University, Im Neuenheimer Feld 205,
69120 Heidelberg, Germany.
Email: \texttt{petra@math.uni-heidelberg.de}.}

\thanks{The authors acknowledge support by the German Research Foundation under Germany’s Excellence Strategy EXC-2181/1 - 390900948 (the Heidelberg STRUCTURES Excellence Cluster). 
}

\begin{document}

\begin{abstract}
 We study the geometric meaning of relative smoothness on dually flat manifolds and its implications for the convergence analysis of mirror descent. Our main result yields one-sided counterparts of the classical Levi--Civita characterizations of Lipschitz Riemannian gradients, each equivalent to relative smoothness.
 We further investigate the stationarity of mirror descent iterates for possibly nonconvex objectives and show that relative smoothness alone does not guarantee convergence of the Riemannian gradient norm to zero. We identify a one-sided condition on the cubic tensor that ensures such convergence and yields an $O(1/N)$ stationarity bound. Examples based on the negative entropy potential illustrate the role of metric variation.
 \end{abstract}
	
	\maketitle
	
	
	
\section{Introduction} \label{sec:introduction}
The classical analysis of gradient descent assumes that the gradient of the objective function $f$ is Lipschitz continuous. Under this assumption, the error in the first-order approximation of $f$ at $x$ admits a quadratic upper bound in the Euclidean displacement $y-x$.
Relative smoothness replaces the squared Euclidean distance by the Bregman divergence associated with a convex potential $\vphi$~\cite{Bauschke:2017aa,Lu2018}. For twice differentiable $f$ and $\vphi$ on an open convex domain $M$, the condition is equivalent to the curvature comparison
\begin{equation*}
    \partial^2 f(x) \preceq L \partial^2 \vphi(x) \qquad \forall x \in M.
\end{equation*}
The choice of $\vphi$ thus allows the upper bound on the curvature of $f$ to vary over the domain. Relative smoothness therefore accommodates objectives whose gradients are not globally Lipschitz continuous. Such objectives arise, for example, in maximum likelihood estimation for Poisson inverse problems~\cite{Bauschke:2017aa,Hurault:23,Daniele2026} and in $D$-optimal experimental design~\cite{Lu2018} 
The resulting bound on the first-order approximation error is fundamental to the convergence analysis of Bregman first-order methods, in particular mirror descent, as it yields a sufficient decrease estimate for a suitable choice of step size.

The same potential also determines the geometry of the mirror descent step. If $\vphi$ is smooth with positive-definite Hessian, then $g=\partial^2\vphi$ defines a Riemannian metric on $M$. The flat connection $\D$ of the primal coordinates $x$ and its dual connection $\Ddual$ with respect to $g$ give $M$ a dually flat structure, with dual affine coordinates $\eta=\partial\vphi(x)$~\cite{amari2016information}. In these coordinates, mirror descent takes the form of natural gradient descent~\cite{RaskuttiMukherjee2015}. Intrinsically, the step follows the $\Ddual$-geodesic whose initial velocity is the negative Riemannian gradient scaled by the step size, so that the dual exponential map serves as a retraction~\cite{RausElshiatyPetra2024}. This interpretation leads us to ask how relative smoothness can be characterized directly in terms of the Riemannian gradient $\ggrad f$ and the dual connection $\Ddual$. Our main result, Theorem~\ref{thm:rel-smth-summary}, shows that the classical Levi--Civita characterizations of Lipschitz Riemannian gradients in terms of Hessian bounds, gradient pairings, and Taylor estimates admit one-sided counterparts on a dually flat manifold, each equivalent to relative smoothness. 
\subsection*{Related Work}
An early form of relative smoothness appeared in \cite{BirnbaumDevanurXiao2011}  in the analysis of Bregman gradient methods for Fisher markets.
Relative smoothness was later systematically developed in \cite{Lu2018,Bauschke:2017aa} as a generalization of Lipschitz gradient smoothness based on a convex reference function.
These works established equivalent convexity, Hessian, and first-order characterizations of relative smoothness, providing a framework for the convergence analysis of Bregman first-order methods beyond classical Lipschitz gradient assumptions. The framework was extended to nonconvex optimization using the notion of $L$-smooth adaptability (L-SMAD) in \cite{BolteVaisbourd2018}, with applications to quadratic inverse problems. 

Relative smoothness is closely connected to the analysis of Bregman proximal gradient methods and mirror descent. The geometric interpretation of mirror descent was investigated in \cite{RaskuttiMukherjee2015}, where mirror descent was identified with natural gradient descent on a dual Riemannian manifold.
We adopt this geometric viewpoint to investigate relative smoothness on dually flat manifolds.

\subsection*{Contribution and Organisation}
To the best of our knowledge, neither a systematic dual-connection characterization of relative smoothness nor a comparison with the classical Levi--Civita characterizations of Lipschitz Riemannian gradients has previously been given. The main contributions of this work are as follows.
\begin{enumerate}
    \item \textit{Geometric characterization of relative smoothness.}
    \begin{enumerate}
        \item \textit{Equivalent formulations.}
        Our central result, Theorem~\ref{thm:rel-smth-summary}, characterizes $L$-relative smoothness by the one-sided bound $g(\Ddual_v\ggrad f,v)\leq L\,g(v,v)$ on the dual covariant derivative of the Riemannian gradient. We establish equivalent formulations in terms of a primal Hessian bound, a gradient pairing controlled by the symmetric Bregman divergence, and a Taylor estimate controlled by the canonical divergence; see Section~\ref{sec:proof-main-result}.
        
        \item \textit{Comparison with Levi--Civita transport.}
        These formulations parallel the classical Levi--Civita characterizations of Lipschitz Riemannian gradients~\cite[Section~10.4]{Boumal:2023}. Dual parallel transport, however, is generally not isometric. We show why the classical Lipschitz estimate for the transported gradient difference does not carry over directly and quantify the effect of metric distortion on this estimate; see Section~\ref{sec:transport}.
    \end{enumerate}
    \item \textit{Riemannian stationarity of mirror descent.}
    \begin{enumerate}
        \item \textit{Stationarity in the averaged metric.}
        We express the classical sufficient-decrease estimate~\cite{Bauschke:2017aa,Lu2018} through the symmetric Bregman divergence and show that the squared gradient norm in the metric averaged along the mirror step tends to zero; see Section~\ref{sec:md-averaged-metric}.
        \item \textit{Sufficient condition for stationarity in the Riemannian norm.}
        We show that the Riemannian gradient norm at the iterates tends to zero under a one-sided bound on the cubic tensor. This bound controls how much the metric can shrink in the gradient direction along each mirror step; see Section~\ref{sec:md-cubic-variation}.
    \end{enumerate}
\end{enumerate}

The paper is organized as follows. Section~\ref{sec:preliminaries} recalls the classical
Levi--Civita characterizations and introduces the dually flat geometric setting and associated divergences. Section~\ref{sec:geom-char-of-rel-smth} establishes
the geometric characterization of relative smoothness and examines the role of dual parallel transport. Section~\ref{sec:MD-as-RGD} investigates the implications
for mirror descent, including Riemannian stationarity under cubic control and examples involving the negative entropy potential.
\subsection*{Notation}
We use $\partial$ for Euclidean derivatives in affine coordinates and $\nabla$ for covariant derivatives. We write $\Exp$ for the exponential map, $\Log$ for its inverse, wherever defined, and $P_{s\to t}$ for parallel transport along a specified curve from parameter $s$ to $t$. Superscripts distinguish the connections used in our analysis, with $\nabla^{\LC}$ denoting the Levi--Civita connection.

For tangent vectors $u,v\in T_xM$, we write $g_x(u,v)$ and $\|u\|_{g(x)}=\sqrt{g_x(u,u)}$. For vector fields $X, Y$ we denote the corresponding pointwise quantities by  $g(X,Y)$ and $\|X\|_g$. We use $\langle\cdot,\cdot\rangle$ for Euclidean pairings in affine coordinates and $\diag(v)$ for the diagonal matrix with vector $v$ on its diagonal.

\section{Geometric Setting}\label{sec:preliminaries}
\subsection{Lipschitz continuous gradients} \label{sec:LC-lipschitz}
For $f \in C^2(M)$, we recall the classical Levi--Civita characterizations of Lipschitz continuity of its Riemannian gradient, specializing the general vector-field formulation in~\cite[Section~10.4]{Boumal:2023}. We also include an equivalent pairing characterization.
Let $L \geq 0$. For each $(x,v)$ in the domain of $\Exp^{\LC}$, let $\gamma(t) := \Exp^{\LC}_x(tv)$, $t \in [0,1]$. We denote by $P^{\LC}_{s\to t}$ the parallel transport along $\gamma$ from $\gamma(s)$ to $\gamma(t)$.
Then the following conditions are equivalent:
\begin{enumerate}
    \item the transport form \cite[Definition~10.44 and Proposition~10.45]{Boumal:2023}: $\ggrad f$ is $L$-Lipschitz continuous if, for all $(x,v)$ in the domain of $\Exp^{\LC}$,
    \begin{equation}\label{eq:LC-Lipschitz-cond}
        \|P^{\LC}_{1 \to 0} \ggrad f(\gamma(1)) - \ggrad f(x) \|_{g(x)} \leq L \|v\|_{g(x)};
    \end{equation}
    \item the norm form \cite[Proposition~10.46]{Boumal:2023} 
    \begin{equation}\label{eq:LC-norm-form}
        \| \nabla_v^{\LC} \ggrad f \|_g \leq L \|v\|_g \qquad \forall v \in TM ;
    \end{equation}
    \item the Hessian form: $\|\Hess^{\LC} f(x)\|_{\op} \leq L$ for all $x \in M$~\cite[Corollary~10.47]{Boumal:2023}, or, since the Hessian is self-adjoint \cite[Proposition~5.15]{Boumal:2023} and the operator norm of a self-adjoint map is its largest absolute eigenvalue \cite[Section~3.1]{Boumal:2023},
    \begin{equation}\label{eq:LC-Hess-bound}
        |\Hess^{\LC} f(v,v)| \leq L\, g(v,v) \qquad \forall v \in TM;
    \end{equation}
    \item the Taylor bound \cite[Corollary~10.54 and Exercise~10.58]{Boumal:2023}: for all $(x,v)$ in the domain of $\Exp^{\LC}$,
    \begin{equation}\label{eq:LC-taylor}
        \bigl|f(\Exp^{\LC}_x(v)) - f(x) - df_x(v)\bigr| \leq \frac{L}{2}\|v\|^2_{g(x)} ;
    \end{equation}
    \item the pairing form: for all $(x,v)$ in the domain of $\Exp^{\LC}$,
    \begin{equation}\label{eq:LC-pairing}
        \bigl| g_{x}\bigl(P^{\LC}_{1 \to 0}\ggrad f(\gamma(1)) - \ggrad f(x),\ v\bigr)\bigr|
        \leq L \|v\|^2_{g(x)} .
    \end{equation}
    Its equivalence with~\eqref{eq:LC-Hess-bound}, and thus with the other characterizations, is proved in Appendix~\ref{app:LC-pairing}.
    \end{enumerate}
Dropping the absolute values in~\eqref{eq:LC-Hess-bound}, \eqref{eq:LC-taylor} and~\eqref{eq:LC-pairing} gives three mutually equivalent statements, all equivalent to the one-sided bound $\Hess^{\LC} f \leq L g$ but no longer to~\eqref{eq:LC-Lipschitz-cond} and~\eqref{eq:LC-norm-form}. 
If $\Hess^{\LC} f$ is positive semidefinite, this distinction disappears.

On dually flat manifolds, the additional affine structure allows us to revisit these characterizations using the dual connection instead, as we shall see in Section~\ref{sec:geom-char-of-rel-smth}.
\subsection{Dually flat manifolds} \label{sec:dually-flat}

Throughout, let $(M,g,\D,\Ddual)$ be a smooth dually flat manifold admitting a global $\nabla$-affine coordinate chart $x$ with open and convex image in $\R^n$. We identify $M$ with this image and use $x$ to denote primal coordinates. In these coordinates, the metric is represented by
\begin{equation*}
    g(x)=\partial^2\vphi(x)\succ0
\end{equation*}
for a smooth Hessian potential $\vphi$. The gradient map $\partial\vphi$ then defines a global $\nabla^\ast$-affine coordinate chart
\begin{equation*}
    \eta=\partial\vphi(x),\qquad M^\ast=\partial\vphi(M),
\end{equation*}
and we assume that $M^\ast$ is also convex.
Writing $\vphi^\ast$ for the Legendre conjugate of $\vphi$, the inverse coordinate map is $x=\partial\vphi^\ast(\eta)$, and the metric in dual coordinates is represented by
\begin{equation*}
    g^\ast(\eta)
    := \partial^2\vphi^\ast(\eta)
    = g\bigl(\partial\vphi^\ast(\eta)\bigr)^{-1}.
\end{equation*}

Thus any two points of $M$ are joined by unique primal and dual geodesic segments, represented by straight segments in their respective affine coordinates.

These assumptions include the smooth, nondegenerate Legendre setting underlying mirror descent. Indeed, if $\vphi:\R^n \to(-\infty,+\infty]$ is a proper lower semicontinuous Legendre function whose restriction to $M=\intr \dom \vphi$ is smooth with positive-definite Hessian, then
\begin{equation*}
    \partial \vphi: M \longrightarrow M^\ast = \intr \dom \vphi^\ast
\end{equation*}
is a diffeomorphism with inverse $\partial \vphi^\ast$ \cite[Theorem~26.5]{Rockafellar:1970} by the positive definiteness of the Hessian and the inverse function theorem. In particular, both coordinate domains are open and convex.

We recall the basic notions of dualistic differential geometry, following~\cite[Sections~6.1 and~6.5]{amari2016information}. 
Let $\VF{M}$ be the space of vector fields on $M$, and $X,Y,Z \in \mathfrak{X}(M)$. 
\begin{itemize}
    \item The \emph{Riemannian gradient} $\ggrad f$ is the vector field on $M$ satisfying
    \begin{equation} \label{eq:grad-def}
        g(\ggrad f, Y) = df(Y).
    \end{equation}
    \item The \emph{Hessian of $f$ with respect to an affine connection $\D$} is the covariant derivative of $df$,
    \begin{equation}\label{eq:def-Riem-Hessian}
        \Hess^{\D}f(X,Y)
        := (\D_X df)(Y)
        = X(df(Y))-df(\D_X Y).
    \end{equation}
    In primal affine coordinates, the Christoffel symbols of $\D$ vanish, so the Hessian is represented by the matrix of second partial derivatives,
    \begin{equation}\label{eq:flat-Hess-is-Eucl}
        \Hess^{\D}f = \partial^2 f.
    \end{equation}
    \item The \emph{covariant derivative of the metric} is
    \begin{equation}\label{eq:nabla-g}
        \nabla g(X,Y,Z) = Xg(Y,Z) - g(\nabla_X Y, Z) - g(Y, \nabla_X Z)
    \end{equation}
    \item The covariant derivatives $\D, \Ddual$ are \emph{dual} if and only if
    \begin{equation} \label{eq:duality}
        Xg(Y,Z) = g(\Ddual_X Y, Z) + g(Y, \D_X Z).
    \end{equation}
    Consequently, the parallel transports along any curve $\gamma$ are mutually adjoint, i.e.\
    \begin{equation} \label{eq:adjoint-transports}
        g_{\gamma(0)}\bigl(P^{\Ddual}_{1\to 0}u,\ w\bigr) = g_{\gamma(1)}\bigl(u,\ P^{\D}_{0\to 1}w\bigr), \qquad u \in T_{\gamma(1)}M,\ w \in T_{\gamma(0)}M .
    \end{equation}
    \item The \emph{cubic tensor} of the conjugate pair $(\D,\Ddual)$ is
    \begin{equation*}
        C(X,Y,Z) := g(\Ddual_X Y-\D_X Y,Z).
    \end{equation*}
    By the duality relation~\eqref{eq:duality} and the expression for the covariant derivative of the metric~\eqref{eq:nabla-g},
    \begin{equation}\label{eq:C-is-nabla-g}
        C = \D g = -\Ddual g.
    \end{equation}
    Since $\D$ and $\Ddual$ are torsion-free and dual, $C$ is totally symmetric. In primal affine coordinates, $C=\partial^3\vphi$. Thus $C=0$ if and only if $g$ is parallel for both connections, equivalently their parallel transports coincide and are isometries. This occurs precisely when $\D=\Ddual=\nabla^{\LC}$.
\end{itemize}
\subsection{Divergences as geodesic energies} \label{sec:segments}
Along a $\nabla^{\LC}$-geodesic $\gamma$ with $\gamma(0)=x$ and $\dot \gamma(0) = v$, the size of a step is unambiguous: since $\gamma$ has constant speed, its energy is the squared norm of its initial velocity, i.e.\ $\int_0^1 \| \dot \gamma(t) \|^2_g dt = \|v\|^2_{g(x)}$ with weighted energy $\int_0^1 (1-t) \|\dot \gamma(t)\|^2_g dt = \tfrac{1}{2}\|v\|^2_{g(x)}$. These are the right-hand sides of~\eqref{eq:LC-pairing} and~\eqref{eq:LC-taylor}. On a dually flat manifold, primal and dual geodesics generally do not have constant speed and the norm of the initial velocity no longer measures the step. We show that the symmetric divergence is the energy shared by the primal and the dual geodesic, and that it is a squared norm of the displacement in a metric averaged along the segment.
The canonical divergence of $(M, g, \D, \Ddual)$ is the Bregman divergence associated with $\vphi$,
\begin{equation*}
    D(x,y) := \vphi(x) - \vphi(y) - \la \partial\vphi(y),\ x - y\ra,
\end{equation*}
and we denote its symmetrization by
\begin{equation*} 
    \Dsym(x,y) := D(x,y) + D(y,x) = \la y - x,\ \partial\vphi(y) - \partial\vphi(x)\ra .
\end{equation*}
Let $x,y \in M$ have dual coordinates $\eta_x = \partial \vphi(x)$ and $\eta_y = \partial \vphi(y)$. The primal geodesic $\gamma^{\D}$ and the dual geodesic $\gamma^{\Ddual}$ from $x$ to $y$, affinely parametrized on $[0,1]$, are segments in their respective affine coordinates,
\begin{equation*}
    \gamma^{\D}(t)=x+t(y-x),
    \qquad
    \partial\varphi\bigl(\gamma^{\Ddual}(t)\bigr)
    =\eta_x+t(\eta_y-\eta_x),
\end{equation*}
We identify $\gamma^{\D}$ with its primal coordinate representation and $\gamma^{\Ddual}$ with its dual coordinate representation $\partial\vphi\circ\gamma^{\Ddual}$ whenever coordinate expressions are used.
Their velocities in these coordinates are the primal and dual displacements
\begin{equation}\label{eq:primal-dual-displacements}
    u := y - x, \qquad a := \eta_y - \eta_x .
\end{equation}
In the other chart, the two geodesics are $\partial \vphi(\gamma^{\D}(t))$ and $\partial \vphi^\ast(\gamma^{\Ddual}(t))$, with velocities
\begin{equation} \label{eq:segment-velocities}
    \frac{d}{dt}\partial \vphi\bigl(\gamma^{\D}(t)\bigr) = g\bigl(\gamma^{\D}(t)\bigr)\,u,
    \qquad
    \frac{d}{dt}\partial \vphi^\ast\bigl(\gamma^{\Ddual}(t)\bigr) = g^\ast\bigl(\gamma^{\Ddual}(t)\bigr)\,a .
\end{equation}
The canonical divergence is the weighted energy of the primal geodesic: by~\cite[Proposition~4]{AyAmari:2015} and~\cite[Section~5.3]{AyAmari:2015}, $D(p,q) = \int_0^1 t\,\|\dot\gamma_{p,q}(t)\|^2_g\,dt$ for the $\nabla$-geodesic $\gamma_{p,q}$ from $p$ to $q$. Applied to the reversal $t \mapsto \gamma^{\D}(1-t)$ of $\gamma^{\D}$, from $y$ to $x$, this reads
\begin{equation*}
    D(y,x) = \int_0^1 (1-t)\,\|\dot\gamma^{\D}(t)\|^2_g\,dt .
\end{equation*}
The symmetric divergence is the full energy, and it is shared by both geodesics.

The following identity is known; see \cite[Theorem~3.2]{amari2016information} for the symmetrized KL-divergence and \cite[Proposition~1]{nielsen2022symmetrized} for the general Bregman case.
\begin{lemma}[Symmetric divergence as shared geodesic energy] \label{lem:shared-energy}
    For $x, y \in M$, with $u$ and $a$ as in~\eqref{eq:primal-dual-displacements},
    \begin{equation*}
        \Dsym(x, y)
        = \la a,\, u\ra
        = \int_0^1 \|\dot\gamma^{\D}(t)\|^2_{g}\,dt
        = \int_0^1 \|\dot\gamma^{\Ddual}(t)\|^2_{g}\,dt .
    \end{equation*}
\end{lemma}
A proof is included in Appendix~\ref{app:shared-energy} for completeness.\\

Lemma~\ref{lem:shared-energy} measures the displacement from $x$ to $y$ by the energy accumulated along the connecting geodesic, rather than by a norm taken at $x$. Equivalently, the symmetric divergence is the squared dual displacement in the metric \emph{averaged} along the dual segment.
\begin{lemma}[Averaged metric form of the symmetric divergence] \label{lm:Dsym-is-dual-displacement-in-avM}
    For $x, y \in M$, with $u$ and $a$ as in~\eqref{eq:primal-dual-displacements}, define 
    \begin{equation*}
        \avM(x,y) := \int_0^1 g^\ast\bigl(\gamma^{\Ddual}(t)\bigr)\, dt.
    \end{equation*}
    Then $\avM(x,y)$ is a metric, $\avM(x,y)\,a = u$, and $\Dsym(x,y) = \|a\|^2_{\avM(x,y)}$.
\end{lemma}
\begin{proof}
   $\avM(x,y)$ is symmetric and positive definite, since $g^\ast$ is a positive definite metric along the dual segment. The identity $\avM(x,y)\,a = u$ follows by integrating the second identity in~\eqref{eq:segment-velocities}; pairing it with $a$ and using Lemma~\ref{lem:shared-energy} gives $\Dsym(x,y) = \la a, u\ra = \|a\|^2_{\avM(x,y)}$.
\end{proof}
In particular, $\avM$ satisfies the secant condition for $\vphi^\ast$: it maps the dual displacement to the primal one exactly, where $g^\ast(\eta_x)$ does so only to first order.
The averaged metric $\avM(x,y)$ coincides with the frozen metric $g^\ast(\eta_x)$ when $g^\ast$ is constant along the dual segment. In general it does not: along the dual segment the displacement $a$ is constant, but by~\eqref{eq:moving-speed} its squared length $\|a\|^2_{g^\ast(\gamma^{\Ddual}(t))}$ equals $\|a\|^2_{g^\ast(\eta_x)}$ at $t = 0$ but may vary with $t$.

The rate of that change is given by the negative of the cubic tensor evaluated on the dual-geodesic velocity,
\begin{equation} \label{eq:speed-drift}
    \frac{d}{dt}\,\|a\|^2_{g^\ast(\gamma^{\Ddual}(t))}
    = \partial^3\vphi^\ast\bigl(\gamma^{\Ddual}(t)\bigr)[a, a, a].
\end{equation}

\section{Geometric Characterization of Relative Smoothness}\label{sec:geom-char-of-rel-smth}
We say that a smooth function $f: M \to \R$ is relatively $L$-smooth ($L\geq 0$) to $\vphi$ iff
\begin{equation}\label{eq:rel-smth-def}
    \partial^2 f(x) \leq L \partial^2 \vphi(x) \quad \text{for all } x \in M
\end{equation}
or equivalently for all $x, y \in M$:
\begin{subequations}
    \begin{align}
        &D_f(x,y) \leq L D_\vphi(x,y) \tag{Bregman inequality} \label{eq:rel-smth-bregman} \\
        &f(x) \leq f(y) + \la \partial f(y), x - y \ra + L D_{\vphi}(x,y). \tag{Taylor expansion} \label{eq:rel-smth-taylor}
    \end{align}
\end{subequations}
Here
\begin{equation*}
    D_f(x,y) := f(x)-f(y)-\la \partial f(y),x-y\ra
\end{equation*}
denotes the Bregman remainder of $f$, which may be negative if $f$ is nonconvex. For convex $f$, this is the Bregman divergence of $f$. In particular, $D_\vphi=D$ is the canonical divergence of Section~\ref{sec:segments}. These coordinate characterizations are classical; see \cite[Proposition~1]{Bauschke:2017aa} and \cite[Proposition~1.1]{Lu2018}. Whenever convexity of $f$ is assumed, it refers to convexity in primal affine coordinates, equivalently $\Hess^{\nabla}f=\partial^2 f\succeq 0$.
We characterize relative smoothness in terms of the Riemannian gradient and the dual connection $\Ddual$. Our main result is the following.
\begin{theorem}[Dual-connection characterization of relative smoothness]
\label{thm:rel-smth-summary}
    Let $f: M \to \R$ be smooth and consider the $\nabla^\ast$-counterparts of the Levi--Civita characterizations in Section~\ref{sec:LC-lipschitz}: for all $x, y \in M$ and $(x,v) \in TM$,
    \begin{enumerate}
        \item[(N)] the norm form $\|\nabla^\ast_v \ggrad f(x)\|_{g(x)} \leq L\|v\|_{g(x)}$;
        \item[(H)] the Hessian bound $|g_x(\nabla^\ast_v \ggrad f(x), v)| \leq L\|v\|^2_{g(x)}$;
        \item[(P)] the pairing form $|g_x(P^{\nabla^\ast}_{1\to 0}\ggrad f(y) - \ggrad f(x),\ y - x)| \leq L\,\Dsym(x,y)$;
        \item[(T)] the Taylor form $|f(y) - f(x) - df_x(y-x)| \leq L\,D(y,x)$.
\end{enumerate}
   Then $f$ is $L$-smooth relative to $\vphi$ if and only if any one of (H), (P) and (T) holds without absolute values.
    If $f$ is convex, this is also equivalent to any of (N), (H), (P) and (T) as stated.
\end{theorem}
\begin{figure}
    \centering
    \includegraphics[width=1.0\linewidth]{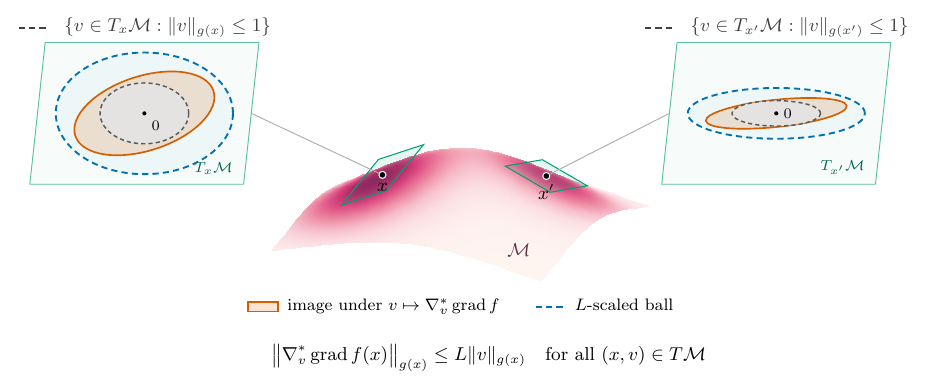}
    \caption{The norm form~(N) of Theorem~\ref{thm:rel-smth-summary}. It states that, at each point, the image of the metric unit ball (grey dashed boundary) under the linear map $v\mapsto\Ddual_v\ggrad f$ is contained in the $L$-scaled ball (blue dashed boundary). The image is shown in orange. For convex $f$, this condition is equivalent to $L$-relative smoothness.}
    \label{fig:norm-form}
\end{figure}
Without convexity, the four conditions (N), (H), (P), and (T) as stated remain equivalent, but characterize the two-sided bound $-Lg\preceq\partial^2 f\preceq Lg$, known as $L$-smooth adaptability (L-SMAD); see \cite[Lemma~2.1]{BolteVaisbourd2018}. For convex $f$, this condition is equivalent to relative smoothness.
Section~\ref{sec:proof-main-result} proves the theorem, first for the pointwise forms (H) and (N) (Section~\ref{sec:hess-norm-forms}), then for the integrated forms (P) and (T) (Section~\ref{sec:geometric-forms})  through a finite-difference characterization for (P). Section~\ref{sec:transport} explains why the classical transport characterization does not extend directly and quantifies the resulting metric distortion.
\begin{remark}[Local extensions]
    The statements of Theorem~\ref{thm:rel-smth-summary} extend locally to dually flat manifolds equipped with local affine charts instead of a global one. On a primal affine coordinate neighborhood $U$ with local Hessian potential $\vphi_U$,
    \begin{equation*}
        \Hess^{\D} f \leq Lg
        \quad \Longleftrightarrow \quad
        \partial^2 f \leq L \partial^2 \vphi_U
        \quad \text{on }U.
    \end{equation*}
    If $U$ has convex images in both primal and dual affine coordinates, the connecting primal and dual geodesic segments remain in $U$, and all equivalences in Theorem~\ref{thm:rel-smth-summary} hold on $U$, with convexity of $f$ understood as convexity of $f|_U$. We work in the global affine-coordinate setting to match the mirror-descent framework considered in Section~\ref{sec:MD-as-RGD}; see also~\cite{Nielsen2020}.
\end{remark}
\subsection{\texorpdfstring{Proof of Theorem~\ref{thm:rel-smth-summary}}{Proof of the main theorem}} \label{sec:proof-main-result}
\subsubsection{Pointwise forms: (H) and (N)} \label{sec:hess-norm-forms}
The next two propositions are closely related. Proposition~\ref{prop:dual-coderivative-of-grad} is an identity that holds for every smooth $f$: it says \emph{which} geometric object the Euclidean Hessian $\partial^2 f$ represents, namely $\nabla^\ast\ggrad f$. Proposition~\ref{prop:rel-smth-is-bounded-Dgrad} then translates the relative smoothness \emph{inequality}~\eqref{eq:rel-smth-def} through this identity.
\begin{proposition} \label{prop:dual-coderivative-of-grad}
On a dually flat Hessian manifold, for every vector field $X$,
\begin{equation*}
    \nabla_X^* \ggrad_g f = g^{-1}\bigl(\Hess^\nabla f(X,\cdot)\bigr).
\end{equation*}
In $\nabla$-affine coordinates, this is represented by
\begin{equation} \label{eq:primal-coord-of-ast-gradf}
    \nabla^* \ggrad_g f = g^{-1}\partial^2 f = (\partial^2 \vphi)^{-1}\partial^2 f.
\end{equation}
\end{proposition}
An equivalent coordinate representation is given in~\cite[Theorem~IV.1]{zhou2025dualriemannian}. A coordinate-free proof based on the duality of the connections is given in Appendix~\ref{app:dual-coderivative-of-grad}. 
\begin{proposition}[Relative Smoothness as an upper bound on the dual covariant derivative]\label{prop:rel-smth-is-bounded-Dgrad}
    A smooth function $f: M \to \R$ on a dually flat Hessian manifold is $L$-smooth relative to the Hessian potential $\vphi$ in primal affine coordinates if and only if
    \begin{equation} \label{eq:quadratic-norm-bound-L-smth}
        g_x(\nabla_v^\ast \ggrad f(x), v) \leq L g_x(v,v) \quad \forall (x,v) \in TM.
    \end{equation}
\end{proposition}
\begin{proof}
    By Proposition~\ref{prop:dual-coderivative-of-grad} and $g=\partial^2 \vphi$,
    \begin{equation*}
        g_x(\nabla_v^\ast \ggrad f(x),v)
        = \la v,\partial^2 f(x)v\ra,
        \qquad
        g_x(v,v)=\la v,\partial^2\vphi(x)v\ra.
    \end{equation*}
    Hence,~\eqref{eq:quadratic-norm-bound-L-smth} is equivalent to~\eqref{eq:rel-smth-def}.
\end{proof}
Proposition~\ref{prop:rel-smth-is-bounded-Dgrad} bounds only the quadratic form of $\nabla^\ast\ggrad f$, and only from above, whereas the norm form (N) in Theorem~\ref{thm:rel-smth-summary}  bounds the operator itself. For convex $f$, the two conditions are equivalent.

\begin{corollary}[Convex case: operator norm and directional norm bounds] \label{cor:cnvx-case}
    Assume in addition that $f$ is convex. Then relative smoothness is equivalent to the operator-norm bound
    \begin{equation} \label{eq:operator-norm-bound}
        \| \nabla^\ast \ggrad f(x) \|_{g(x), \op} \leq L \qquad \forall x \in M,
    \end{equation}
    and thus to the directional bound
    \begin{equation} \label{eq:directional-norm-bound}
        \| \nabla_v^\ast \ggrad f(x)\|_{g(x)} \leq L \|v\|_{g(x)} \qquad \forall (x,v)\in TM .
    \end{equation}
\end{corollary}
\begin{proof}
    By Proposition~\ref{prop:rel-smth-is-bounded-Dgrad}, relative $L$-smoothness is equivalent to
    \begin{equation*}
        g_x(\nabla_v^\ast \ggrad f(x), v) \leq L g_x(v,v)
        \qquad \forall (x,v)\in TM.
    \end{equation*}
    Moreover, $\nabla^\ast \ggrad f = g^{-1}\partial^2 f$ is $g$-self-adjoint. Since $f$ is convex, $\partial^2 f$ is positive semidefinite, hence $\nabla^\ast \ggrad f(x)$ is $g(x)$-positive semidefinite for every $x$. Therefore all eigenvalues are nonnegative, and
    \begin{equation}\label{eq:cnvx-EV}
        \| \nabla^\ast \ggrad f(x)\|_{g(x), \op}
        =
        \lambda_{\max}\bigl(\nabla^\ast \ggrad f(x)\bigr).
    \end{equation}
    Without convexity, we are only guaranteed $\| \cdot \|_{g,\op} = \max_i |\lambda_i(\cdot)|$, which may exceed $\lambda_{\max}$; semi-definiteness is what ensures \eqref{eq:cnvx-EV}. It follows that the quadratic-form bound
    \begin{equation*}
        g_x(\nabla_v^\ast \ggrad f(x), v) \leq L g_x(v,v)
        \qquad \forall (x,v)\in TM
    \end{equation*}
    is equivalent to the operator norm bound \eqref{eq:operator-norm-bound}. The equivalence of \eqref{eq:operator-norm-bound} and \eqref{eq:directional-norm-bound} is the definition of the operator norm.
\end{proof}
For convex $f$, $\partial^2 f \succeq 0$, so $g_x(\nabla^\ast_v\ggrad f(x), v) = \la v, \partial^2 f(x)\, v\ra \geq 0$ by Proposition~\ref{prop:dual-coderivative-of-grad}, and Proposition~\ref{prop:rel-smth-is-bounded-Dgrad} holds with absolute values, i.e.\ in the two-sided form~\eqref{eq:LC-Hess-bound} with $\nabla^\ast$ in place of $\nabla^{\LC}$.
\subsubsection{Integrated forms: (P) and (T)} \label{sec:geometric-forms}
We first express relative smoothness as a second-derivative bound along primal geodesics. Integrating this bound once yields the finite-difference characterization and hence the pairing form~(P); integrating it twice yields the Taylor form~(T).
\begin{lemma}[Second-derivative characterization of relative smoothness] \label{lem:rel-smth-second-deriv}
    Let $(M, g, \nabla, \nabla^\ast)$ be a dually flat Hessian manifold with canonical divergence $D$. Let $f: M \to \R$ be a smooth function. For any primal geodesic $\gamma(t) = \Exp_x^\nabla(tv)$ define
    \begin{equation} \label{eq:geodesic-functions-for-second-variation}
        h(t) := f(\gamma(t)), \quad \quad \quad q(t):= D(\gamma(t), \gamma(0)).
    \end{equation}
    Along every such geodesic,
    \begin{equation*}
        h''(t)
        = \Hess^\nabla f(\dot\gamma(t),\dot\gamma(t)),
        \qquad
        q''(t) = \|\dot\gamma(t)\|_g^2,
        \qquad
        q(0)=q'(0)=0.
    \end{equation*}
    Consequently, $f$ is $L$-smooth relative to $\vphi$ if and only if
    \begin{equation*}
        h''(t) \leq L q''(t)
    \end{equation*}
    for every primal geodesic $\gamma$ and every $t$ for which $\gamma(t)$ is defined in $M$.
\end{lemma}
The proof is given in Appendix~\ref{app:relative-smoothness}
\begin{proposition}[Finite-difference form of relative smoothness] \label{prop:finite-diff-form}
    Let $f: M \to \R$ be smooth. Then $f$ is $L$-smooth relative to $\vphi$ if and only if for every primal geodesic $\gamma(t) = \Exp_x^\nabla(tv)$,
    \begin{equation} \label{eq:finite-diff-form}
        df_{\gamma(1)}(\dot\gamma(1)) - df_{\gamma(0)}(\dot\gamma(0)) \leq L \Dsym(\gamma(1), \gamma(0)).
    \end{equation}
\end{proposition}
\begin{proof}
Assume first that $f$ is $L$-smooth relative to $\vphi$. With $h$ and $q$ as in~\eqref{eq:geodesic-functions-for-second-variation}, $h'(t) = df_{\gamma(t)}(\dot\gamma(t))$, so integrating $h'' \leq L q''$ from Lemma~\ref{lem:rel-smth-second-deriv} gives
\begin{equation*}
    df_{\gamma(1)}(\dot\gamma(1)) - df_{\gamma(0)}(\dot\gamma(0))
    = \int_0^1 h''(t)\,dt
    \;\leq\; L\int_0^1 q''(t)\,dt
    = L\int_0^1 \|\dot\gamma(t)\|^2_g\,dt,
\end{equation*}
and by Lemma~\ref{lem:shared-energy} the remaining integral equals $\Dsym(\gamma(1), \gamma(0))$.

Conversely, fix a primal geodesic $\gamma$ and $t$. For $s > 0$ small enough that $\gamma(t+s)$ is defined, consider the shifted and rescaled geodesic
\begin{equation*}
    r \mapsto \Exp^\nabla_{\gamma(t)}\bigl(r s\,\dot\gamma(t)\bigr) = \gamma(t + rs), \qquad r \in [0,1],
\end{equation*}
from $\gamma(t)$ to $\gamma(t+s)$, with velocities $s\,\dot\gamma(t)$ at $r = 0$ and $s\,\dot\gamma(t+s)$ at $r = 1$. Applying~\eqref{eq:finite-diff-form} to it, the left-hand side is $s\bigl(h'(t+s) - h'(t)\bigr)$.
By Lemma~\ref{lem:rel-smth-second-deriv} 
\begin{equation*}
    q'(t+s)-q'(t) = \int_t^{t+s}\|\dot\gamma(\tau)\|_g^2\,d\tau.
\end{equation*}
Thus, the right-hand side is
\begin{equation*}
    L\Dsym\bigl(\gamma(t+s), \gamma(t)\bigr) 
    = L \int_0^1 \|s\,\dot\gamma(t+rs)\|^2_g\,dr
    = L s \int_t^{t+s} \|\dot\gamma(\tau)\|^2_g\,d\tau
    = L s\bigl(q'(t+s) - q'(t)\bigr),
\end{equation*}
where we used the change of variables 
$\tau = t + rs$.
Hence
\begin{equation*}
    \frac{h'(t+s) - h'(t)}{s} \leq L\,\frac{q'(t+s) - q'(t)}{s}.
\end{equation*}
Passing to the limit as $s \to 0$ gives $h''(t) \leq L q''(t)$ for every primal geodesic $\gamma$ and every $t$, and Lemma~\ref{lem:rel-smth-second-deriv} concludes that $f$ is $L$-smooth relative to $\vphi$.
\end{proof}
\begin{corollary}[Pairing form] \label{cor:LC-pairing-one-sided}
    $f$ is $L$-smooth relative to $\vphi$ if and only if, for all $x, y \in M$, 
    with $\gamma^{\D}$ the primal geodesic from $x$ to $y$
    and $v := \dot\gamma^{\D}(0)$,
    \begin{equation} \label{eq:pairing-form}
        g_x\bigl(P^{\Ddual}_{1\to 0}\ggrad f(y) - \ggrad f(x),\ v\bigr) \leq L \Dsym(x,y),
        \qquad
        \Dsym(x,y) = \int_0^1 \|\dot\gamma^{\D}(t)\|^2_g\,dt .
    \end{equation}
\end{corollary}
\begin{proof}
    The velocity of the $\nabla$-geodesic $\gamma^{\D}$ is $\nabla$-parallel, so $\dot\gamma^{\D}(1) = P^{\nabla}_{0\to 1}v$, and by the adjointness~\eqref{eq:adjoint-transports} of the dual transports,
    \begin{equation*}
        g_x\bigl(P^{\nabla^\ast}_{1\to 0}\ggrad f(y),\ v\bigr) = g_y\bigl(\ggrad f(y),\ \dot\gamma^{\D}(1)\bigr) = df_{\gamma(1)}(\dot \gamma^{\D}(1)).
    \end{equation*}
    Hence the claim is Proposition~\ref{prop:finite-diff-form}.
\end{proof}
\begin{corollary}[Taylor expansion form] \label{cor:taylor-expansion-form}
    Let $f: M \to \R$ be smooth. Then $f$ is $L$-smooth relative to $\vphi$ if and only if, for all $(x,v) \in TM$ with $\Exp^\nabla_x(v) \in M$,
    \begin{equation*}
        f(\Exp^\nabla_x(v)) \leq f(x) + df_x(v) + L\,D(\Exp^\nabla_x(v), x),
        \qquad
        D(\Exp^\nabla_x(v), x) = \int_0^1 (1-t)\,\|\dot\gamma(t)\|^2_g\,dt,
    \end{equation*}
    where $\gamma(t) = \Exp^\nabla_x(tv)$.
\end{corollary}
\begin{proof}
    Assume first that $f$ is $L$-smooth relative to $\vphi$. Recall the definitions of $h(t)$ and $q(t)$ from \eqref{eq:geodesic-functions-for-second-variation}, and the one-dimensional integral remainder formula
    \begin{equation*}
        h(1)-h(0)-h'(0)=\int_0^1 (1-t)h''(t)dt, \qquad q(1)-q(0)-q'(0)=\int_0^1 (1-t)q''(t)dt.
    \end{equation*}
    Since $q(0)=0$ and $q'(0)=0$, Lemma~\ref{lem:rel-smth-second-deriv} yields
    \begin{equation*}
        h(1)-h(0)-h'(0)\le L q(1).        
    \end{equation*}
    Substituting the definitions of $h$ and $q$ gives
    \begin{equation*}
        f(\Exp^\nabla_x(v)) \leq f(x) + df_x(v) + L\,D(\Exp^\nabla_x(v), x).
    \end{equation*}
    Conversely, fix a primal geodesic $\gamma$ and $t$. Define the shifted geodesic
    \begin{equation*}
        \gamma_t(s):= \Exp_{\gamma(t)}^\nabla(s \dot \gamma(t)) = \gamma(t+s),
    \end{equation*}
    where the last equality holds since $\gamma$ is an affinely parametrized $\nabla$-geodesic. Set
    \begin{equation*}
        h_t(s):=f(\gamma_t(s)), \qquad q_t(s):=D(\gamma_t(s),\gamma_t(0)).
    \end{equation*}
    Applying the assumed inequality at the base point $\gamma(t)$ with tangent vector $s\,\dot\gamma(t)$ gives
    \begin{equation*}
        h_t(s)\leq h_t(0)+s\,h_t'(0)+Lq_t(s).
    \end{equation*}
    Hence
    \begin{equation*}
        \frac{h_t(s)-h_t(0)-s\,h_t'(0)}{s^2} \leq L \frac{q_t(s)}{s^2}.
    \end{equation*}
    Passing to the limit as $s \to 0$ and using $q_t(0)=q_t'(0)=0$, we get
    \begin{equation*}
        \frac{1}{2} h_t''(0) \leq \frac{L}{2} q_t''(0).
    \end{equation*}
    That is,
    \begin{equation*}
        h''(t) \leq L q''(t) \qquad \forall \gamma, \forall t.
    \end{equation*}
    Applying Lemma~\ref{lem:rel-smth-second-deriv} concludes that $f$ is $L$-smooth relative to $\vphi$.
\end{proof}
\begin{proof}[Proof of Theorem~\ref{thm:rel-smth-summary}]
    By Proposition~\ref{prop:rel-smth-is-bounded-Dgrad}, Corollary~\ref{cor:LC-pairing-one-sided}, and Corollary~\ref{cor:taylor-expansion-form} (using $\Exp_x^\nabla(y-x)=y$), relative smoothness is equivalent to each of (H), (P), and (T) without absolute values.

    If $f$ is convex, then $\Hess^\nabla f=\partial^2 f\succeq 0$, so $h''\geq 0$ in Lemma~\ref{lem:rel-smth-second-deriv}. Consequently, the quadratic form in (H), the pairing in (P), and the Taylor remainder in (T) are nonnegative. Thus each of these conditions is equivalent to its version with absolute values. The equivalence with (N) follows from Corollary~\ref{cor:cnvx-case}.
\end{proof}
\begin{remark}[Intrinsic forms of the Euclidean characterizations] \label{rem:euclidean-forms}
    Each of the equivalent Euclidean forms of the definition has an intrinsic counterpart on the dually flat manifold: the Hessian bound~\eqref{eq:rel-smth-def} becomes the quadratic-form bound on $\nabla^\ast\ggrad f$ (Proposition~\ref{prop:rel-smth-is-bounded-Dgrad}), and the Taylor expansion~\eqref{eq:rel-smth-taylor} becomes its geodesic version, Corollary~\ref{cor:taylor-expansion-form}. The latter also recovers the Bregman inequality~\eqref{eq:rel-smth-bregman}: in primal affine coordinates $\Log^\nabla_x y = y-x$, hence
    \begin{equation*}
        f(y) - f(x) - df_x(\Log^\nabla_x y) = D_f(y,x).
    \end{equation*}
    Thus, Corollary~\ref{cor:taylor-expansion-form} reads
    \begin{equation*}
        D_f(y,x) \leq L D_\vphi(y,x) = L D(y,x).
    \end{equation*}
    
\end{remark}
\subsection{The Transport Form}\label{sec:transport}
Theorem~\ref{thm:rel-smth-summary} contains a $\Ddual$-counterpart of every Levi--Civita characterization in Section~\ref{sec:LC-lipschitz} except the transport form~\eqref{eq:LC-Lipschitz-cond}. This section explains why the corresponding transport condition is strictly stronger than relative smoothness even for convex functions. 
We proceed in the following steps:
\begin{enumerate}
    \item we state a literal $\Ddual$-transcription of the transport form~\eqref{eq:LC-Lipschitz-cond} and show that it implies relative smoothness (Proposition~\ref{prop:transport-implies-rel-smooth});
    \item we give a convex, relatively smooth function that violates it;
    \item we quantify the failure (Proposition~\ref{prop:transport-up-to-kappa}): relative smoothness implies the literal transport form up to a factor $\kappa$ measuring the distortion of the metric along the segment. This factor is unbounded in the example, whereas $\Dsym$ absorbs the distortion and needs no such factor. 
\end{enumerate}
\subsubsection{The $\Ddual$-transcription of~\eqref{eq:LC-Lipschitz-cond}}
We first record what the Riemannian gradient looks like in dual coordinates.
\begin{lemma}[Dual-coordinate representation of $\nabla^\ast \ggrad_g f$] \label{cor:dual-primal-operator-norms-agree}
    Let $\eta = \partial \vphi(x)$ be the $\nabla^\ast$-affine coordinates and write $\ggrad_g f = \sum_{i=1}^n G^i(\eta)\,\partial_{\eta_i}$. Then, 
    \begin{equation} \label{eq:G-is-euclidean-gradient}
        G(\eta) = \partial_x f\bigl(\partial\vphi^\ast(\eta)\bigr),
    \end{equation}
    i.e.\ in dual coordinates the Riemannian gradient has as components the Euclidean gradient of $f$ in primal coordinates. Moreover, $DG(\eta)$ represents
    $\nabla^\ast\ggrad_g f$ in these coordinates, and
    \begin{equation} \label{eq:op-norms-agree}
        \|DG(\eta)\|_{g^\ast(\eta),\op} = \|\nabla^\ast \ggrad_g f(x)\|_{g(x),\op} .
    \end{equation}
\end{lemma}
\begin{proof}
    For~\eqref{eq:G-is-euclidean-gradient}, vector components transform as $V^i_\eta = \sum_j \frac{\partial \eta_i}{\partial x_j} V^j_x$, and $\frac{\partial\eta_i}{\partial x_j} = \partial_i\partial_j\vphi = g_{ij}$ since $\eta = \partial\vphi(x)$. With $\ggrad_g f = g^{-1}\partial_x f$ in primal components,
    \begin{equation*}
        G^i(\eta) = \sum_{j,k} g_{ij}\,(g^{-1})^{jk}\,\partial_{x_k} f = \partial_{x_i} f .
    \end{equation*}
    Since $\eta$ are $\nabla^\ast$-affine coordinates, the Christoffel symbols of $\nabla^\ast$ vanish in these coordinates, so for every vector field $X = \sum_j X^j \partial_{\eta_j}$,
    \begin{equation*}
        \nabla_X^\ast \ggrad_g f = \sum_{i,j} X^j\, \partial_{\eta_j} G^i(\eta)\,\partial_{\eta_i},
    \end{equation*}
    and the matrix of the $(1,1)$-tensor $\nabla^\ast \ggrad_g f$ in the $\eta$-coordinates is $DG(\eta)$. Finally, $g^\ast(\eta) = \partial^2\vphi^\ast(\eta)$ is the coordinate representation of the metric in these coordinates, so the operator norm of $DG(\eta)$ with respect to $g^\ast(\eta)$ is the pointwise operator norm of $\nabla^\ast\ggrad_g f$ with respect to $g(x)$, which is~\eqref{eq:op-norms-agree}.
\end{proof}
Since $\eta$ are $\nabla^\ast$-affine, $\nabla^\ast$-parallel transport preserves $\eta$-components. For $x,y \in M$, by~\eqref{eq:G-is-euclidean-gradient}, $G(\eta_y) - G(\eta_x)$ therefore represents $P^{\nabla^\ast}_{1\to 0}\,\ggrad_g f(y) - \ggrad_g f(x)$, with transport along the dual geodesic $\gamma^{\Ddual}$, whose velocity in dual coordinates is the displacement $a = \eta_y - \eta_x$. The
$\nabla^\ast$-transcription of the transport form~\eqref{eq:LC-Lipschitz-cond} is thus
\begin{equation} \label{eq:transport-dual}
    \|G(\eta_y) - G(\eta_x)\|_{g^\ast(\eta_x)} \leq L\,\|a\|_{g^\ast(\eta_x)},
    \qquad \forall x,y \in M.
\end{equation}
\begin{proposition} \label{prop:transport-implies-rel-smooth}
    If~\eqref{eq:transport-dual} holds, then $f$ is $L$-smooth relative to $\vphi$.
\end{proposition}
\begin{proof}
    We fix $x$. Since $\eta$ are affine coordinates, for $v \in \R^n$ 
    and $|t|$ small, $\eta_x + tv = \eta_y$ for some $y \in M$, with $a = tv$. Then
    \begin{equation*}
        \|G(\eta_x +tv) - G(\eta_x) \|_{g^\ast(\eta_x)} \leq L |t| \,\|v\|_{g^\ast(\eta_x)}.
    \end{equation*}
    Dividing by $|t|$ turns the left-hand side into the norm of a finite difference. Since $G$ is differentiable, letting $t \to 0$ gives
    \begin{equation*}
        \frac{G(\eta_x +tv) - G(\eta_x)}{t} \to DG(\eta_x)v,
    \end{equation*}
    and hence
    \begin{equation*}
        \| DG(\eta_x) v\|_{g^\ast(\eta_x)} \leq L \|v\|_{g^\ast(\eta_x)}
        \qquad \forall v\in \R^n,
    \end{equation*}
    so $\|DG(\eta_x)\|_{g^\ast(\eta_x),\op}\leq L$.
    By Lemma~\ref{cor:dual-primal-operator-norms-agree}, this is $\|\nabla^\ast \ggrad_g f(x)\|_{g(x),\op}\leq L$, and by Cauchy--Schwarz,
    \begin{equation*}
        g_x(\nabla_v^\ast \ggrad_g f(x),v) \leq \|\nabla_v^\ast \ggrad_g f(x)\|_{g(x)}\,\|v\|_{g(x)} \leq L\|v\|_{g(x)}^2
    \end{equation*}
    for all $(x,v)\in TM$. Proposition~\ref{prop:rel-smth-is-bounded-Dgrad} now implies that $f$ is $L$-smooth relative to $\vphi$.
\end{proof}
\subsubsection{A counterexample}
The converse of Proposition~\ref{prop:transport-implies-rel-smooth} is false, even for convex functions.
\begin{example} \label{ex:transport-counterexample}
    Let $\vphi(x) = \sum_{i=1}^2 ( x_i \ln x_i - x_i )$ on $M = \R^2_{++}$, so that
    \begin{equation*}
        \eta = \ln x, \qquad \vphi^\ast(\eta) = e^{\eta_1} + e^{\eta_2},
    \end{equation*}
    with $M^\ast = \R^2$ and $g^\ast(\eta) = \diag(x)$. Furthermore, let $f(x) := s \ln s$ with $s := x_1 + x_2$. Both $M$ and $M^\ast$ are convex.
    The function $f$ is convex and $1$-smooth relative to $\vphi$: this follows by a Cauchy--Schwarz argument using $\partial^2 f(x) = \tfrac{1}{s}\,\mathbf{1}\mathbf{1}^\top$:
    \begin{equation*}
        \frac{(v_1+v_2)^2}{x_1+x_2} \leq \frac{v_1^2}{x_1} + \frac{v_2^2}{x_2}, \quad \forall v \in \R^2 \cong T_xM,
    \end{equation*}
    i.e.\ $0 \preceq \partial^2 f \preceq g$ on all of $M$. Now, by~\eqref{eq:G-is-euclidean-gradient},
    \begin{equation*}
        G(\eta) = \partial_x f = (1+ \ln s)\,\mathbf{1},
        \qquad s = e^{\eta_1} + e^{\eta_2} .
    \end{equation*}
    For $K \geq 1$, let $x, y \in M$ with dual coordinates $\eta_x = (0,-K)$ and $\eta_y = (0,0)$, i.e.\ $x = (1, e^{-K})$, $y = (1,1)$ and $a = (0,K)$. Then
    \begin{equation*}
        \bigl\|G(\eta_y) - G(\eta_x)\bigr\|_{g^\ast(\eta_x)}
        = \sqrt{1+e^{-K}}\, \ln\left(\frac{2}{1+e^{-K}}\right),
    \end{equation*}
    and $\|a\|_{g^\ast(\eta_x)} = K e^{-K/2}$. Therefore, using $\sqrt{1+e^{-K}} \geq 1$ and $1 + e^{-K} \leq 1 + e^{-1}$,
    \begin{equation*}
        \frac{\|G(\eta_y) - G(\eta_x)\|_{g^\ast(\eta_x)}}{\|a\|_{g^\ast(\eta_x)}}
        \geq \frac{e^{K/2}}{K} \ln\left(\frac{2}{1+e^{-1}}\right)  
        \;\xrightarrow[K\to\infty]{}\; \infty .
    \end{equation*}
    Hence~\eqref{eq:transport-dual} fails for every $L$, and so does the converse of Proposition~\ref{prop:transport-implies-rel-smooth}.
\end{example}
The violation depends on which point serves as the base point: interchanging $x$ and $y$, the ratio becomes $\sqrt{2}\,\ln\frac{2}{1+e^{-K}}/K \to 0$. The blow-up in Example~\ref{ex:transport-counterexample} comes from the base point $x = (1, e^{-K})$ approaching the boundary of $M$, where $g^\ast(\eta_x)$ becomes singular. In the Levi--Civita setting this asymmetry cannot occur, since parallel transport is an isometry and~\eqref{eq:LC-Lipschitz-cond} is symmetric in its endpoints. Here both sides of~\eqref{eq:transport-dual} are measured in the frozen metric $g^\ast(\eta_x)$, and the failure is governed by how much $g^\ast$ changes between the two endpoints, which Proposition~\ref{prop:transport-up-to-kappa} quantifies.
\begin{remark}[The importance of the right metric]\label{rem:transport-vs-pairing}
    For the Levi--Civita connection, the transport form~\eqref{eq:LC-Lipschitz-cond} and pairing form~\eqref{eq:LC-pairing} are equivalent. For $\Ddual$, both have a pairing version, and these differ only in the metric. Let $x, y \in M$. By Cauchy--Schwarz, the literal transport form~\eqref{eq:transport-dual} implies
    \begin{equation*}
        \la G(\eta_y) - G(\eta_x),\ g^\ast(\eta_x)\,a\ra \;\leq\; L\,\|a\|^2_{g^\ast(\eta_x)} .
    \end{equation*}
    In the pairing form~\eqref{eq:pairing-form}, the transported gradient difference has $\eta$-components $G(\eta_y) - G(\eta_x)$, as noted above, while $v = \dot\gamma^{\D}(0)$ has primal components $u$, hence $\eta$-components $g(x)\,u$. Since $g^\ast(\eta_x) = g(x)^{-1}$ represents the metric in dual coordinates,
    \begin{equation*}
        g_x\bigl(P^{\Ddual}_{1\to 0}\ggrad f(y) - \ggrad f(x),\ v\bigr)
        = \la G(\eta_y) - G(\eta_x),\ g^\ast(\eta_x)\, g(x)\, u\ra
        = \la G(\eta_y) - G(\eta_x),\ u\ra ,
    \end{equation*}
    and the secant identity $u = \avM(x,y)\,a$ of Lemma~\ref{lm:Dsym-is-dual-displacement-in-avM} turns~\eqref{eq:pairing-form} into
    \begin{equation*}
        \la G(\eta_y) - G(\eta_x),\ \avM(x,y)\,a\ra
        \leq L\|a\|^2_{\avM(x,y)}
        = L\Dsym(x,y) .
    \end{equation*}
    The pairing form uses the averaged metric $\avM$, which samples $g^\ast$ at every point of the dual segment, whereas the literal transport form evaluates $g^\ast$ at $\eta_x$ only. If $g^\ast$ is constant along the dual segment, then $\avM(x,y)=g^\ast(\eta_x)$ and the two pairing inequalities coincide. This mirrors the Levi--Civita case, where parallel transport preserves the metric, so the moving and frozen metrics agree after transport. In general, the averaged pairing accounts for the variation of the metric tensor along the dual segment, whereas the literal transport pairing uses the metric frozen at $\eta_x$.
\end{remark}
\subsubsection{Quantifying the failure}
Remark~\ref{rem:transport-vs-pairing} locates the obstruction; the following proposition makes it quantitative. It is not a restatement of the pairing form (P): it keeps the frozen metric of the literal transport form~\eqref{eq:transport-dual} and asks which constant relative smoothness yields there. The answer is $L$ times a factor $\kappa$ that records the distortion between the frozen and the moving metric along the segment.
\begin{proposition}\label{prop:transport-up-to-kappa}
    Let $f$ be convex and $L$-smooth relative to $\vphi$. For $x,y\in M$, set $\eta_x=\partial\vphi(x)$, $\eta_y=\partial\vphi(y)$, and let $\gamma^{\Ddual}(t) = (1-t)\eta_x + t \eta_y$ be the dual geodesic from $x$ to $y$. Then
    \begin{equation*}
        \|G(\eta_y) - G(\eta_x)\|_{g^\ast(\eta_x)} \leq L\, \kappa(\eta_x, \eta_y)\, \|\eta_y-\eta_x\|_{g^\ast(\eta_x)},
        \qquad \kappa(\eta_x,\eta_y) := \int_0^1 \rho(t)\sigma(t)\, dt \geq 1,
    \end{equation*}
    with $\rho(t) = \|g^\ast(\eta_x)^{1/2}g^\ast(\gamma^{\Ddual}(t))^{-1/2}\|_{\op}$ and $\sigma(t) = \|g^\ast(\gamma^{\Ddual}(t))^{1/2}g^\ast(\eta_x)^{-1/2}\|_{\op}$ the forward and backward distortions between the frozen and the moving metric along the segment.
\end{proposition}
\begin{proof}
    Set $a:=\eta_y-\eta_x$. By convexity and Lemma~\ref{cor:dual-primal-operator-norms-agree}, relative smoothness implies
    \begin{equation*}
    \|DG(\gamma^{\Ddual}(t) )\|_{g^\ast(\gamma^{\Ddual}(t)), \op} \leq L,\qquad t\in[0,1].
    \end{equation*}
    The statement then follows by the fundamental theorem of calculus and the definitions of $\rho$ and $\sigma$,
       \begin{align*}
           \|G(\eta_y) - G(\eta_x)\|_{g^\ast(\eta_x)}
           &\leq \int_0^1 \|DG(\gamma^{\Ddual}(t)) \, a\|_{g^\ast(\eta_x)}\, dt \\
           &\leq \int_0^1 \rho(t)\, \|DG(\gamma^{\Ddual}(t))\, a\|_{g^\ast(\gamma^{\Ddual}(t))}\, dt \\
           &\leq L \int_0^1 \rho(t)\, \| a \|_{g^\ast(\gamma^{\Ddual}(t))}\, dt \\
           &\leq L \int_0^1 \rho(t) \sigma(t)\, dt \; \| a \|_{g^\ast(\eta_x)} .
       \end{align*}
\end{proof}
The two factors account for the two changes of metric in the proof: $\rho$ transfers the integrand from the moving metric to the frozen one, and $\sigma$ transfers the step from the frozen metric to the moving one. The naive bound~\eqref{eq:transport-dual} is recovered when $\kappa \equiv 1$, which happens if and only if  $g^\ast(\gamma^{\Ddual}(t))$ is a scalar multiple of $g^\ast(\eta_x)$ for all $t$; in particular when $\nabla^\ast g = 0$ along the segment.

In Example~\ref{ex:transport-counterexample}, $\gamma^{\Ddual}(t) = (0, -(1-t)K)$ in dual coordinates, so  $g^\ast(\eta_x)^{1/2} g^\ast(\gamma^{\Ddual}(t))^{-1/2} = \diag(1, e^{-tK/2})$, hence $\rho \equiv 1$, $\sigma(t) = e^{tK/2}$ and
\begin{equation*}
    \kappa(\eta_x, \eta_y) = \int_0^1 e^{tK/2}\,dt = \frac{2\,(e^{K/2} - 1)}{K}.
\end{equation*}
Since $\kappa$ is unbounded as $K\to\infty$, the estimate does not yield a uniform transport bound of the form~\eqref{eq:transport-dual}.
The pairing form (P) needs no such factor: measuring the step by $\Dsym = \|a\|^2_{\avM}$ incorporates the variation of the metric along the entire dual segment, without requiring an additional distortion factor.
\section{Mirror Descent as Riemannian Gradient Descent}\label{sec:MD-as-RGD}
We revisit mirror descent from the viewpoint of dually flat information geometry \cite{RaskuttiMukherjee2015,RausElshiatyPetra2024}. In this setting, the geometric characterization of relative smoothness in Theorem~\ref{thm:rel-smth-summary} is particularly natural, since mirror descent steps follow geodesics of the dual connection $\Ddual$.

We first recall the classical sufficient-decrease result and express it in terms of the symmetric Bregman divergence, or equivalently, the metric averaged along the step.
We then investigate what relative smoothness implies for stationarity. While it guarantees decay of the averaged Bregman energy, this alone does not ensure convergence of the Riemannian gradient norm to zero. To this end, we introduce an additional geometric condition controlling the variation of the metric along mirror descent steps.
\subsection{Descent in the averaged metric} \label{sec:md-averaged-metric}
Let $\alpha > 0$. The mirror step is defined by
\begin{equation*}
    x_{k+1} = T_\alpha(x_k) := \Exp^{\nabla^\ast}_{x_k}\bigl(-\alpha \ggrad f(x_k)\bigr),
\end{equation*}
or, equivalently, $\eta_{k+1} = \eta_k - \alpha\, G(\eta_k)$ in dual coordinates, since by~\eqref{eq:G-is-euclidean-gradient} the $\eta$-components of $\ggrad f$ are $G(\eta) = \partial_x f(x)$. We assume throughout that the dual steps are admissible, i.e.\ $\eta_k-\alpha G(\eta_k)\in M^\ast$ for all $k$. Consequently, $x_{k+1}\in M$. Mirror descent is thus Riemannian gradient descent along the dual geodesic. $T_\alpha(x_k)$ uniquely minimizes $y \mapsto df_{x_k} (\Log^{\nabla}_{x_k} y) + \frac{1}{\alpha} D(y,x_k)$. This is directly evident when writing the map in dual coordinates: the objective is $\la \partial_x f(x_k), y - x_k \ra + \frac{1}{\alpha} D(y,x_k)$, whose stationarity condition is $\partial_x f(x_k) + \frac{1}{\alpha} (\eta_y - \eta_k) = 0$.
\subsubsection*{Sufficient descent}
At $y = T_\alpha(x_k)$, using Lemma~\ref{lem:shared-energy}, we have
\begin{equation*}
    \alpha\, df_{x_k}(\Log^{\nabla}_{x_k} x_{k+1})
    = \alpha \la G(\eta_k),\, x_{k+1} - x_k\ra
    = -\Dsym(x_k, x_{k+1}) .
\end{equation*}
Note that this identity holds independently of relative smoothness. Applying the descent form of relative smoothness (Corollary~\ref{cor:taylor-expansion-form}) at $x_k$ with the primal velocity $\Log^{\nabla}_{x_k} x_{k+1}$ and substituting this identity yields
\begin{equation*}
    f(x_{k+1}) \leq f(x_k) - \frac{1}{\alpha} \Dsym(x_k, x_{k+1}) + L\, D(x_{k+1}, x_k).
\end{equation*}
Since $D(x_{k+1}, x_k) \leq \Dsym(x_k, x_{k+1})$, this gives the sufficient-descent estimate
\begin{equation} \label{eq:md-sufficient-descent}
    f(x_k) - f(x_{k+1}) \geq \Bigl(\tfrac{1}{\alpha} - L\Bigr)\Dsym(x_k, x_{k+1}),
\end{equation}
so for any step size $\alpha>0$ satisfying $\alpha L<1$, the objective values $f(x_k)$ are monotonically nonincreasing.
\subsubsection*{The decrease in the averaged metric.} By Lemma~\ref{lm:Dsym-is-dual-displacement-in-avM}
\begin{equation*}
    \Dsym(x_k, x_{k+1}) = \alpha^2\, \|G(\eta_k)\|^2_{\avM_k},
    \qquad \avM_k := \avM(x_k, x_{k+1}) .
\end{equation*}
The sufficient decrease~\eqref{eq:md-sufficient-descent} is therefore a squared gradient norm, measured in the metric averaged along the step:
\begin{equation} \label{eq:md-decrease-H}
    f(x_k) - f(x_{k+1}) \;\geq\; \alpha(1 - \alpha L)\, \|G(\eta_k)\|^2_{\avM_k} .
\end{equation}
This motivates the stationarity measure
\begin{equation*}
    \mathcal{S}_\alpha(x) := \frac{1}{\alpha^2}\Dsym\bigl(x, T_\alpha(x)\bigr)
    = \|G(\eta)\|^2_{\avM(x,\,T_\alpha(x))} ,
\end{equation*}
which vanishes exactly at stationarity points, i.e.\ $\ggrad f(x) = 0$. Moreover, as $\alpha \to 0$, it reduces to the squared Riemannian gradient norm,
\begin{equation*}
    \mathcal{S}_\alpha(x)
    \longrightarrow
    \|G(\eta)\|^2_{g^\ast(\eta)}
    =
    \|\ggrad f(x)\|^2_{g(x)}.
\end{equation*}
The sufficient-descent estimate follows the classical relative-smoothness argument \cite{Bauschke:2017aa,Lu2018}. The following proposition uses this estimate to establish a stationary bound in the averaged metric.
\begin{proposition}[Convergence of mirror descent]\label{prop:md-convergence}
Let $f$ be $L$-smooth relative to $\vphi$ and bounded below, and suppose that the mirror steps are admissible. Let $\alpha > 0$, $\alpha L < 1$. Then
\begin{equation*}
    \sum_{k=0}^{\infty} \|G(\eta_k)\|^2_{\avM_k}
    = \sum_{k=0}^{\infty} \mathcal{S}_\alpha(x_k)
    \;\leq\;
    \frac{f(x_0) - \inf f}{\alpha(1 - \alpha L)}
    \;<\; \infty ,
\end{equation*}
so $\mathcal{S}_\alpha(x_k) \to 0$ and every accumulation point of $(x_k)_k$ in $M$ is stationary.
\end{proposition}
\begin{proof}
    Summing up~\eqref{eq:md-decrease-H} yields the claimed bound. If $x_{k_j}\to\bar x\in M$, then $\Dsym(x_{k_j},x_{k_j+1})\to0$. Since $\bar x$ is an interior point, the positive-definite Hessian of $\vphi$ is uniformly bounded below in a neighborhood of $\bar x$. Integrating the Hessian along the segment and using convexity of $\vphi$ therefore yields $x_{k_j+1}-x_{k_j}\to0$. The mirror-descent update then gives $G(\eta_{k_j})\to0$, thus $\ggrad f(\bar x)=0$.
\end{proof}
\subsubsection*{Convergence in the Riemannian norm.} Since $g^\ast(\eta)$ represents the metric in dual coordinates, $\|G(\eta)\|^2_{g^\ast(\eta)} = \|\ggrad f(x)\|^2_{g(x)}$. The averaged and the Riemannian norm of the gradient differ by the ratio
\begin{equation*}
    \Lambda_k := \frac{\|G(\eta_k)\|^2_{\avM_k}}{\|G(\eta_k)\|^2_{g^\ast(\eta_k)}},
\end{equation*}
where we set $\Lambda_k=1$ if $G(\eta_k)=0$.
The ratio equals one when the metric is constant along the step. A sufficient condition for obtaining convergence of the Riemannian gradient norm from Proposition~\ref{prop:md-convergence} is that $\Lambda_k$ remains bounded away from zero. Geometrically, this means that the metric does not shrink too much along the step in the direction of the gradient.
\begin{corollary} \label{cor:md-riemannian}
Under the assumptions of Proposition~\ref{prop:md-convergence}, suppose $\Lambda_k \geq \lambda > 0$ for all $k$. Then
\begin{equation*}
    \sum_{k=0}^{\infty} \|\ggrad f(x_k)\|^2_{g(x_k)}
    \;\leq\; \frac{f(x_0) - \inf f}{\lambda\,\alpha(1 - \alpha L)} ,
\end{equation*}
and in particular $\|\ggrad f(x_k)\|_{g(x_k)} \to 0$.
\end{corollary}
\begin{proof}
    The result follows from $\mathcal S_\alpha(x_k) = \Lambda_k\|\ggrad f(x_k)\|^2_{g(x_k)}$ and Proposition~\ref{prop:md-convergence}.
\end{proof}
\subsection{Cubic variation and Riemannian stationarity}\label{sec:md-cubic-variation}
 We now give a geometric sufficient condition for the
uniform lower bound $\Lambda_k \geq \lambda > 0$. Along the mirror step, written in dual coordinates as $\gamma^{\Ddual}(t)=\eta_k-\alpha tG(\eta_k)$, define
\begin{equation*}
    q_k(t):=\|G(\eta_k)\|^2_{g^\ast(\gamma^{\Ddual}(t))}.
\end{equation*}
This is the squared norm of the initial gradient under $\Ddual$-parallel transport. Its initial value and average satisfy
\begin{equation*}
    q_k(0)=\|\ggrad f(x_k)\|^2_{g(x_k)},
    \qquad
    \int_0^1q_k(t)\,dt=\mathcal S_\alpha(x_k).
\end{equation*}
Thus, whenever $q_k(0) \neq 0$,
\begin{equation*}
    \Lambda_k=\frac{\int_0^1q_k(t)\,dt}{q_k(0)}.
\end{equation*}
Since dual parallel transport preserves the components $G(\eta_k)$, only the metric varies in $q_k$.
By~\eqref{eq:speed-drift},
\begin{equation}\label{eq:md-cubic-drift}
    q_k'(t)
    =-\alpha\,\partial^3\vphi^\ast
      \bigl(\gamma^{\Ddual}(t)\bigr)
      [G(\eta_k),G(\eta_k),G(\eta_k)].
\end{equation}
In dual affine coordinates, $\partial^3\vphi^\ast$ represents the negative of the cubic tensor $C=\nabla g=-\nabla^\ast g$, see \eqref{eq:C-is-nabla-g}. The following theorem shows that a one-sided bound on this term in the mirror-step direction limits the decrease of $q_k$ and yields a uniform positive lower bound on $\Lambda_k$.
\begin{theorem}[Riemannian stationarity under cubic control] \label{thm:md-cubic-stationarity}
    Under the assumptions of Proposition~\ref{prop:md-convergence}, suppose there exists a constant $K\geq 0$ such that along every mirror descent step, 
    \begin{equation}\label{eq:md-cubic-condition}
        \partial^3\vphi^\ast(\gamma^{\Ddual}(t))
        [G(\eta_k),G(\eta_k),G(\eta_k)]
        \leq K\,q_k(t), \qquad \text{for all} \ t\in[0,1] .
    \end{equation}
    Then
    \begin{equation*}
        \Lambda_k \geq \lambda_{\alpha,K}
        := \frac{1-e^{-\alpha K}}{\alpha K}
        \qquad \text{for all } k,
    \end{equation*}
    where $\lambda_{\alpha,0}=1$ by continuity, so that $\lambda_{\alpha,K}>0$ for all $K$. Consequently, {$\|\ggrad f(x_k)\|_{g(x_k)}\to0$} and
    \begin{equation*}
        \min_{0\leq k<N}\|\ggrad f(x_k)\|^2_{g(x_k)}
        = O(1/N).
    \end{equation*}
\end{theorem}
\begin{proof}
    By~\eqref{eq:md-cubic-drift}, the cubic bound implies
    $q_k'(t)\geq-\alpha K q_k(t)$. Gr\"onwall's inequality
    therefore gives $q_k(t)\geq e^{-\alpha Kt}q_k(0)$.
    For $q_k(0)>0$, it follows that
    \begin{equation*}
        \Lambda_k
        = \int_0^1 \frac{q_k(t)}{q_k(0)}\,dt
        \geq \int_0^1 e^{-\alpha Kt}\,dt
        = \lambda_{\alpha,K}.
    \end{equation*}
    If $q_k(0)=0$, the bound holds by the convention $\Lambda_k=1$.
    The remaining statements follow from Corollary~\ref{cor:md-riemannian} and Proposition~\ref{prop:md-convergence}.
\end{proof}
\subsection{Negative entropy: uniform metric compatibility and its failure} \label{sec:md-negative-entropy}
For negative entropy, the cubic condition~\eqref{eq:md-cubic-condition} can be checked directly in primal coordinates. We provide two examples to illustrate the role of this condition as the iterates approach the solution on the boundary of the domain. The second shows that relative smoothness alone does not ensure stationarity in the Riemannian norm.

Consider the negative entropy potential on $M = \R^n_{++}$
\begin{equation*}
    \vphi(x)=\sum_{i=1}^n(x_i\ln x_i-x_i). 
\end{equation*}
Its dual coordinates, dual potential and dual metric are given by
\begin{equation*}
    \eta=\ln x,\qquad
    \vphi^\ast(\eta)=\sum_{i=1}^n e^{\eta_i},\qquad
    g^\ast(\eta)=\diag(e^\eta).
\end{equation*}
The mirror descent update reads
\begin{equation*}
    x_{k+1}=x_{k}e^{-\alpha \partial_x f(x_k)},
\end{equation*}
where exponential and multiplication are understood componentwise.
\begin{corollary}[Negative entropy] \label{cor:md-entropy}
    Under the assumptions of Proposition~\ref{prop:md-convergence}, suppose there exists $K\geq0$ such that
    \begin{equation}\label{eq:md-entropy-bound}
        \|(\partial_x f(x_k))_+\|_\infty \leq K
        \qquad \text{for all } k\geq0,
    \end{equation}
    where $(\cdot)_+$ denotes the positive part of a vector. Then the cubic condition~\eqref{eq:md-cubic-condition} holds with the same constant $K$, and the conclusions of Theorem~\ref{thm:md-cubic-stationarity} follow.
\end{corollary}
Recall from~\eqref{eq:G-is-euclidean-gradient} that the dual-coordinate components of the Riemannian gradient are given by the Euclidean gradient in primal coordinates, $G(\eta_k) = \partial_x f(x_k)$.
\begin{proof}
    Since $\diag(G(\eta_k))\preceq KI$ and $g^\ast(\eta)$ is diagonal and positive, we have
    \begin{equation*}
        \partial^3\vphi^\ast(\eta) [G(\eta_k),G(\eta_k),G(\eta_k)]
        = G(\eta_k)^\top g^\ast(\eta)\diag(G(\eta_k))G(\eta_k)
        \leq K\,G(\eta_k)^\top g^\ast(\eta)G(\eta_k).
    \end{equation*}
    Along the dual segment of the $k$-th mirror step, the right-hand side is $Kq_k(t)$, which gives \eqref{eq:md-cubic-condition}.
\end{proof}
We illustrate the condition in Corollary~\ref{cor:md-entropy}
and its failure through two one-dimensional examples.
\subsubsection*{Uniform metric compatibility}
For $n=1$, consider $f(x)=x+1-e^{-x}$. Since $f''<0$, relative smoothness holds with $L=0$. Moreover, $1<f'(x)<2$, so~\eqref{eq:md-entropy-bound} holds with $K = 2$. For any $\alpha > 0$, the mirror updates are $x_{k+1}=x_k e^{-\alpha f'(x_k)}\leq x_k e^{-\alpha}$, hence $x_k\to0$ and 
$f'(x_k) \to 2$. Integrating along the dual segment yields
\begin{equation*}
    \Lambda_k
    =\frac{1-e^{-\alpha f'(x_k)}}{\alpha f'(x_k)}
    \longrightarrow
    \frac{1-e^{-2\alpha}}{2\alpha}>0.
\end{equation*}
Thus the averaged metric remains bounded below relative to the metric at $x_k$, even as the iterates approach the boundary.
\subsubsection*{Failure of uniform metric compatibility}
For $f(x)=\sqrt{x}$, we again have $f''(x) < 0$ and relative smoothness with $L = 0$. Here, $f'(x)=\tfrac{1}{2\sqrt x}$ is unbounded as $x\to0$. The mirror update gives $x_{k+1}=x_k e^{-\alpha/(2\sqrt{x_k})}$. With sufficient descent, the factor multiplying $x_k$ is upper bounded $e^{-\alpha/(2\sqrt{x_0})}<1$. Hence, $x_k \to 0$ and
\begin{equation*}
    \Lambda_k
    =\frac{1-e^{-\alpha/(2\sqrt{x_k})}}{\alpha/(2\sqrt{x_k})} \longrightarrow 0,
\end{equation*}
while
\begin{equation*}
    \|\ggrad f(x_k)\|^2_{g(x_k)}=\frac{1}{4}, \qquad \mathcal S_\alpha(x_k) = \frac{\Lambda_k}{4} \longrightarrow 0.
\end{equation*}
Thus, relative smoothness alone does not guarantee convergence of the Riemannina gradient norm to zero, since the averaged metric may fail to remain uniformly comparable to the metric at the starting point $x_k$.
\section{Conclusion}

We have characterized relative smoothness on dually flat manifolds through the dual covariant derivative of the Riemannian gradient. The equivalent Hessian, gradient pairing, and Bregman Taylor formulations parallel the classical Levi--Civita characterizations, while the failure of dual parallel transport to preserve the metric explains why the classical Lipschitz gradient estimate does not carry over directly.

For mirror descent, we distinguish stationarity in the metric averaged along each step from stationarity in the Riemannian metric at the iterates. A one-sided condition on the cubic tensor connects these notions, and negative entropy examples illustrate the role of this condition.

The geometric characterization provides a basis for further analysis of mirror descent and for developing extensions that exploit the dual affine structure and account for metric variation.

\section*{Acknowledgments}
\paragraph{\textbf{Use of AI tools}} The authors used ChatGPT (Sol 5.6) and Claude (Opus 5.5) for language editing (rephrasing, grammar, and synonyms). Furthermore, it was used to generate initial draft of Figure~\ref{fig:norm-form} and to enhance the Example~\ref{ex:transport-counterexample} that were both subsequently revised substantially. The authors reviewed and verified all such outputs and take full responsibility for the content of this paper.
\appendix
\section{Proofs of auxiliary results} \label{app:auxiliary-proofs}
\subsection{Levi--Civita pairing characterization} \label{app:LC-pairing}
\begin{proof}[Proof: Equivalence of the Hessian and pairing forms]
    Fix $(x,v)$ in the domain of $\Exp^{\LC}$, and set $\gamma(t) := \Exp^{\LC}_x(tv)$ and $h := f\circ\gamma$, for $t$ for which $\gamma(t)$ is defined. Since $\dot\gamma$ is parallel along $\gamma$ \cite[(10.25)]{Boumal:2023}, we have
    \begin{equation*}
        h''(t) = \Hess^{\LC} f(\dot\gamma(t),\dot\gamma(t)).
    \end{equation*}
    Moreover, since $P^{\LC}$ is an isometry \cite[Proposition~10.36]{Boumal:2023}, 
    \begin{equation*}
        g_x\bigl(P^{\LC}_{t\to 0}\ggrad f(\gamma(t)) - \ggrad f(x),\,tv\bigr)
        = t\bigl(h'(t)-h'(0)\bigr).
    \end{equation*}
    For $t>0$, the left-hand side of~\eqref{eq:LC-pairing} applied to $tv$ therefore equals $t\,|h'(t)-h'(0)|$.

    Assume first that~\eqref{eq:LC-Hess-bound} holds. We obtain
    \begin{equation*}
        |h'(1)-h'(0)|
        \leq \int_0^1 |\Hess^{\LC} f(\dot\gamma(t),\dot\gamma(t))|\,dt
        \leq L\int_0^1 \|\dot\gamma(t)\|_g^2\,dt
        = L\|v\|_{g(x)}^2,
    \end{equation*}
    since $\gamma$ has constant speed. This proves~\eqref{eq:LC-pairing}.

    Conversely, assume that~\eqref{eq:LC-pairing} holds. Applying it to $tv$ for sufficiently small $t>0$ gives
    \begin{equation*}
        t\,|h'(t)-h'(0)| \leq Lt^2\|v\|_{g(x)}^2.
    \end{equation*}
    Dividing by $t^2$ and letting $t\to 0$ yields
    \begin{equation*}
        |\Hess^{\LC} f(v,v)| = |h''(0)|
        \leq L\|v\|_{g(x)}^2.
    \end{equation*}
    Since every $v\in T_xM$ defines such a geodesic locally, this proves~\eqref{eq:LC-Hess-bound}.
\end{proof}
\subsection{Symmetric divergence and geodesic energy} \label{app:shared-energy}
\begin{proof}[Proof of Lemma~\ref{lem:shared-energy}]
By~\eqref{eq:segment-velocities}, $\|\dot\gamma^{\D}(t)\|^2_g = \la u,\, g(\gamma^{\D}(t))\,u\ra = \bigl\la u,\, \tfrac{d}{dt}\partial\vphi(\gamma^{\D}(t))\bigr\ra$, which integrates to $\la u, a\ra$. For the dual geodesic, its squared speed in primal coordinates is
    \begin{equation}\label{eq:moving-speed}
        \|\dot\gamma^{\Ddual}(t)\|^2_{g}
        = \bigl\la g^\ast(\gamma^{\Ddual}(t))\,a,\ g\bigl(\partial\vphi^\ast(\gamma^{\Ddual}(t))\bigr)\, g^\ast(\gamma^{\Ddual}(t))\,a \bigr\ra
        = \|a\|^2_{g^\ast(\gamma^{\Ddual}(t))} ,
    \end{equation}
    since $g\bigl(\nabla\vphi^\ast(\gamma^{\Ddual}(t))\bigr)\, g^\ast(\gamma^{\Ddual}(t)) = I$.
    Then, integrating the second identity in~\eqref{eq:segment-velocities} similarly gives $\int_0^1 \|\dot\gamma^{\Ddual}(t)\|^2_g\,dt = \la a, u\ra$. Finally, $\la a, u\ra = D(x,y) + D(y,x) = \Dsym(x,y)$.
\end{proof}
\subsection{Dual covariant derivative of the gradient} \label{app:dual-coderivative-of-grad}
\begin{proof}[Proof of Proposition~\ref{prop:dual-coderivative-of-grad}]
    Applying the duality relation~\eqref{eq:duality} to the pair $(\ggrad_g f, Y)$ gives
    \begin{equation*}
        X g(\ggrad_g f, Y) = g(\nabla_X^* \ggrad_g f, Y) + g(\ggrad_g f, \nabla_X Y).
    \end{equation*}
    Since $g(\ggrad_g f, Y)=df(Y)$ by~\eqref{eq:grad-def}, this becomes
    \begin{equation*}
        X df(Y)=df(\nabla_X Y)+g(Y,\nabla_X^* \ggrad_g f).    
    \end{equation*}
    Therefore
    \begin{equation} \label{eq:primal-Hess-is-dual-coderivative-of-grad}
        g(Y,\nabla_X^* \ggrad_g f)
    = Xdf(Y)-df(\nabla_X Y)
    = \Hess^\nabla f(X,Y).
    \end{equation}
    By~\eqref{eq:flat-Hess-is-Eucl} and $g=\partial^2 \vphi$ in $\nabla$-affine coordinates, this implies~\eqref{eq:primal-coord-of-ast-gradf}.
\end{proof}
\subsection{Second derivative characterisation of relative smoothness} \label{app:relative-smoothness}
\begin{proof}[Proof of
  Lemma~\ref{lem:rel-smth-second-deriv}]
    We first compute the second derivative of $q$. By the geodesic-energy formula for the canonical divergence,
    \begin{equation*}
        D(p,q) = \int_0^1 t \| \dot \gamma_{p,q}(t) \|^2_g dt,
    \end{equation*}
    where $\gamma_{p,q}$ is the $\nabla$-geodesic from $p$ to $q$. This canonical divergence was first introduced in \cite[Proposition~4]{AyAmari:2015} and they show in Section~5.3 that it coincides with the canonical divergence defined in terms of the Bregman divergence of $\vphi$ on dually flat manifolds. Since $s \mapsto \gamma(t(1-s))$ is the $\nabla$-geodesic from $\gamma(t)$ to $\gamma(0)$, with velocity $-t\,\dot\gamma(t(1-s))$, we get
    \begin{equation*}
        q(t) = \int_0^1 st^2 \| \dot \gamma(t(1-s))\|^2_g ds.
    \end{equation*}
    With the change of variables $u=(1-s)t$, so that $ds=-du/t$, this becomes
    \begin{equation*}
        q(t) = \int_0^t (t-u) \| \dot \gamma(u) \|^2_g du.
    \end{equation*}
    Differentiating twice yields
    \begin{equation*}
        q^\prime(t) = \int_0^t \| \dot \gamma(u) \|^2_g du \qquad q^{\prime \prime}(t) = \| \dot \gamma(t) \|^2_g.
    \end{equation*}
    Now consider $h(t) = f(\gamma(t))$. Then,
    \begin{equation*}
        h^\prime(t) = df_{\gamma(t)} (\dot \gamma(t)).
    \end{equation*}
    Differentiating once more along $\gamma$ using the primal connection gives
    \begin{equation*}
        h^{\prime \prime}(t) =(\nabla df)_{\gamma(t)}(\dot\gamma(t),\dot\gamma(t)) +df_{\gamma(t)}\bigl(\nabla_{\dot\gamma(t)}\dot\gamma(t)\bigr).
    \end{equation*}
    Since $\gamma$ is a $\nabla$-geodesic, the second term vanishes. Therefore,
    \begin{equation*}
        h^{\prime \prime}(t) = \Hess^\nabla f(\dot \gamma(t), \dot \gamma(t)) =  g(\dot \gamma(t),\nabla_{\dot \gamma(t)}^* \ggrad_g f).
    \end{equation*}
    where the second equality is~\eqref{eq:primal-Hess-is-dual-coderivative-of-grad}.

    Combining the two computations, we have
    \begin{equation*}
        h''(t)\leq Lq''(t) \quad\Longleftrightarrow\quad g_{\gamma(t)}\bigl(\nabla_{\dot\gamma(t)}^*\ggrad_g f,\dot\gamma(t)\bigr) \leq Lg_{\gamma(t)}(\dot\gamma(t),\dot\gamma(t)).
    \end{equation*}
    Since every pair $(x,v)\in TM$ is realized as
    \begin{equation*}
        x=\gamma(0),\qquad v=\dot\gamma(0)    
    \end{equation*}
    for some primal geodesic $\gamma(t)=\Exp_x^\nabla(tv)$, the latter inequality for all $\gamma$ and all $t$ is equivalent to
    \begin{equation*}
        g_x\bigl(\nabla_v^*\ggrad_g f,v\bigr)\leq Lg_x(v,v) \qquad \forall (x,v)\in TM.
    \end{equation*}
    By Proposition~\ref{prop:rel-smth-is-bounded-Dgrad}, this is exactly the characterization of $L$-smoothness of $f$ relative to $\vphi$. This proves the equivalence.
\end{proof}

\printbibliography

\end{document}